# Group decision making with q-rung orthopair hesitant fuzzy preference relations

Benting Wan[1]*  Jiao Zhang[1]

[1]School of Software and IoT Engineering, Jiangxi University of Finance and Economics, Nanchang 330013, China

**Abstract**：This paper mainly studies group decision making (GDM) problem based on q-rung orthopair hesitant fuzzy preference relations (q-ROHFPRs). First, Based on q-rung orthopair hesitant fuzzy sets (q-ROHFSs), the definitions of q-ROHFPR and additive consistent q-ROHFPR are introduced. However, the ideal state of complete consistency does not exist in real life. In order to better judge whether the decision maker's preference matrix satisfies consistency, the consistency index of q-ROHFPR is used to judge whether the matrix of q-ROHFPR is acceptable. For the q-ROHFPR matrix that does not meet the acceptable consistency, two optimization models are established for deriving the acceptably additive consistent q-ROHFPRs. In order to make the q-ROHFPR matrix of decision makers still satisfy the consistency after aggregation, this paper extends the q-rung orthopair hesitant fuzzy weighted geometric average operator (q-ROHFWGA). At the same time, in order to verify whether decision makers can reach consensus after aggregation, a consensus index based on distance is offered. Based on this consensus index, an optimization model that satisfies consistency and consensus is constructed to solve the priority vector, and develop a consistency and consensus-based approach for dealing with group decision-making (GDM) with q-ROHFPRs. Finally, the case in this paper verifies the validity and accuracy of the group decision-making model, and also verifies that the q-ROHFPR consistency and consensus management model proposed in this paper can solve the q-rung orthopair hesitant fuzzy preference group decision-making problem.

## 1 Introduction

As information becomes more and more complex, and the ambiguity of decision-making becomes more and more obvious, decision makers (DMs) may be very strict with quantifying preferences with real number. At the same time, considering the limitations of fuzzy set (FS) [1] in dealing with the clear

---

value that DMs are unwilling to provide their judgment information. Torra [2] introduced the definition of hesitant fuzzy set (HFS), in which the degree of membership is represented by multiple discrete values instead of a single value, which solves the MODM and MOGDM problems more flexibly [ 3][4][5][6]. For example, Li et al. [7] proposed a method for measuring the distance between hesitant fuzzy elements. Xu and Xia [8, 9] studied the aggregation operators, distance and similarity of HFSs. And Zhu et al. [10] introduced dual hesitant fuzzy set (DHFS) and studied its basic operations and properties. Chen et al. [11] put forward interval-valued hesitant fuzzy set (IVHFS) and basic operations and properties, and applied them to cluster analysis. Liu et al. [12] merged HFS with q-rung orthopair fuzzy sets (q-ROFS), proposed q-rung orthopair hesitant fuzzy sets (q-ROHFS), and studied its properties and distance measurement methods. The q-rung orthopair hesitant fuzzy number (q-ROHFN) can well describe the uncertainty and ambiguity in reality. since Satty proposed the multiplicative preference relation (MPR) to deal with the priority of multiple alternatives[13] , Multi-objective decision-making (MODM) based on preference relations has received extensive attention from researchers and has been widely used in production and life [14] [15]. The judgment matrix composed of the preference relations between the alternatives is convenient for expressing people's judgments on the affairs. In order to express the preference relations between the  alternatives more flexibly, the researchers extended the MPR to fuzzy sets, and hesitant fuzzy preference relations (HFPR) [16], dual hesitant fuzzy preference relations (DHFPR) [17], interval hesitant fuzzy preference relations (IVHFPR) [18] and other fuzzy preference relations based on additive consistency and multiplicative consistency are presented, which are used to deal with MODM problems in more complex environments [19][20][21]. After investigation and analysis, q-ROHFNs have not yet been used to express preference relations and solve MODM and MOGDM problems. This paper will study q-rung orthopair hesitant fuzzy preference relation (q-ROHFPR) and its properties.

When DMs use fuzzy preference relations to express judgment information, it is often difficult to ensure consistency. Therefore, judging whether fuzzy preference relations satisfy consistency has been a hot research topic in recent years. For this reason, researchers have provided a large number of optimization models for HFPR, DHFPR, IVHFPR and other preference relations, which are used to judge whether the matrix given by the DMs meets the consistency. Zhang et al. [22] established a minimum deviation optimization model based on consistent HFPR. The model minimizes the number of adjustments and the number of elements that need to be adjusted; Zhu et al. [23] proposed α-normalization and β-normalization methods; in order to obtain the interval consistency index, Li et al. [24] studied and constructed two models based on optimism; Tang et al. [25] developed a programming model to confirm the complete additive consistency of HFPR; Zhang et al. [26] constructed several goal programming models to verify the acceptable additive consistency of HFPR , and then built a programming model to verify the acceptable multiplicative consistency of HFPR [27]; Meng and An [28] developed several 0-1 mixed programming models to derive priority weights vector from the multiplicative consistency of HFPR. In order to be able to use q-ROHFPR to solve the MODM problem conveniently, we will study the q-ROHFPR consistency judgment model and adjustment model. However, when solving the MOGDM problems, not only the consistency of the decision makers' information must be considered, but also the consensus of group decision-making needs to be considered, so that the opinions reserved between the decision makers can reach agreement and consensus. For this reason, researchers have developed optimization methods for consistency and consensus: Wu et al. [29] presented two heuristic algorithms to improve the consistency and consensus of the AHP framework. Zhang et al. [30] used an iterative algorithm to improve the consistency and consensus of the AHP framework. To solve the consistency and consensus of HMPR, He and Xu [31] established a consensus model of hesitant information on the basis of nonlinear programming models with different preference structures. In order to further solve the MOGDM problem of q-ROHFPRs, this paper will further study the consensus model based on the consistency optimization model. The main innovations of this paper are as follows:

(1) The definition of q-ROHFPR is proposed. Based on q-ROHFSs, the concept of q-ROHFPR is introduced, and its six properties are analyzed. Compared with HFPRs, q-ROHFPRs have more relaxed restrictions, which is convenient for DMs to express their preferences more flexibly.

(2) The additive consistency of q-ROHFPRs and the consistency index (CI) of q-ROHFPRs are defined. Inspired by [32], the additive consistent q-ROHFPR is defined. And inspired by the literature 错误!未找到引用源。, the consistency index (CI) of q-ROHFPR is defined.

(3) The acceptably consistent adjustment optimization model and the acceptably additive consistency priority model of q-ROHFPRs are established. Based on CI, an acceptably consistent adjustment optimization model is designed. And inspired by the literature [34], an additive consistency priority weight model of q-ROHFPRs is proposed which is used to derive the q-ROHF priority weight vector based on the acceptably additive consistent q-ROHFPRs matrix, which is used to solve the MOMD problem.

(4) A consensus adjustment model of q-ROHFPRs is developed. Firstly, q-rung orthopair hesitant fuzzy weighted geometric average operator is extended, and a distance-based group consensus index (GCI) is developed. A goal programming model to derive new q-ROHFPR with acceptable consistency and consensus are constructed with minimal adjustment, and it is used to solve the MOGDM problem of q-ROHFPRs. Case results show that proposed method is effective and feasible.

The remainder of the paper is arranged as: Section 2 describes preliminaries. Section 3 explains the concept of q-ROHFPR and additive consistency of q-ROHFPR. Section 4 elaborates the linear goal programming models of acceptably additive consistent q-ROHFPRs. Section 5 conducts consensus analysis on q-ROHFPR. Section 6 uses a case to verify the feasibility and effectiveness of induced model. Section 7 describes our work and future research.

## 2 Preliminaries

In this section, some concepts about $q-ROHFNs$ are reviewed.

**Definition 1[12]**：Let $X = \{x_1, x_2, ..., x_n\}$ be a $q-ROHFS$ $Q$ on $X$ is shown in formula (1).

$$Q = \{< x_i, h(x_i), g(x_i) > | x_i \in X\}, (q \geq 1) \tag{1}$$

In formula (1), $h(x_i)$ and $g(x_i)$ are called all possible membership functions and non-membership functions of $Q$, respectively, $h(x_i)$ and $g(x_i)$ are two non-empty subsets on [0,1] respectively. For all $x \in U$, there are $\forall x_i \in X$, $\forall u \in h(x_i)$, $\forall v \in g(x_i)$, $u^+ = max_{u \in h(x_i)}\{u\}$, $v^+ = max_{v \in g(x_i)}\{v\}$, but $u^+$ and $v^+$ also need to meet the following: $0 \leq (u^+)^q + (v^+)^q \leq 1$, $0 \leq u \leq 1$, $0 \leq v \leq 1, (q \geq 1)$. In addition, the hesitation function of the element $x_i \in X$ can be expressed as: $\pi_Q(x_i) = \sqrt[q]{1 - u^q - v^q}, (q \geq 1)$.

Usually, we use an ordered pair of membership degree and non-membership degree, namely $Q = < h(x), g(x) >_q$ represent the q-rung hesitant fuzzy number $(q - ROHFN)$.

**Definition 2[12]:** Let $Q = <h, g>$、 $Q_1 = <h_1, g_1>$ and $Q_2 = <h_2, g_2>$ be three $q-ROHFNs$, where $(q \geq 1)$, $\lambda \neq 0$ is any real number, then the basic algorithm of $q - ROHFNs$ can be defined as:

(1) $Q_1 \oplus Q_2 = \bigcup_{u_1 \in h_1, u_2 \in h_2, v_1 \in g_1, v_2 \in g_2} < \left\{ \sqrt[q]{u_1^q + u_2^q - u_1^q u_2^q} \right\}, \{v_1 v_2\} >$;

(2) $Q_1 \otimes Q_2 = \bigcup_{u_1 \in h_1, u_2 \in h_2, v_1 \in g_1, v_2 \in g_2} < \{u_1 u_2\}, \left\{ \sqrt[q]{v_1^q + v_2^q - v_1^q v_2^q} \right\} >$;

(3) $\lambda Q = \bigcup_{u \in h, v \in g} < \left\{ \sqrt[q]{1-(1-u^q)^\lambda} \right\}, \{v^\lambda\} >$;

(4) $Q^\lambda = \bigcup_{u \in h, v \in g} < \{u^\lambda\}, \left\{ \sqrt[q]{1-(1-v^q)^\lambda} \right\} >$.

**Definition 3[12]:** Let $Q = <h, g> = <\{u_1, u_2, \ldots, u_{l(h)}\}, \{v_1, v_2, \ldots, v_{l(h)}\} >, (q \geq 1)$ is a $q-ROHFN$, and its score function and accuracy function are shown in formula (2) and formula (3).

$$S_Q = \frac{1}{l(h)} \sum_{u \in h} u^q - \frac{1}{l(g)} \sum_{v \in g} v^q, \tag{2}$$

$$D_Q = \frac{1}{l(h)} \sum_{u \in h} u^q + \frac{1}{l(g)} \sum_{v \in g} v^q, \tag{3}$$

In formulas (2) and (3), $l(h)$ and $l(g)$ represent the number of elements contained in h and g, respectively.

**Definition 4[12]:** Let $Q_1 = <h_1, g_1>$ and $Q_2 = <h_2, g_2>$ be two $q - ROHFNs$, then the size comparison rules are defined as follows:

(1) If $S_{Q_1} > S_{Q_2}$, then $Q_1 > Q_2$;

(2) If $S_{Q_1} < S_{Q_2}$, then $Q_1 < Q_2$;

(3) If $S_{Q_1} = S_{Q_2}$, then calculate the accuracy of the two $q - ROHFNs$ respectively, and compare the accuracy values. If $D_{Q_1} > D_{Q_2}$, then $Q_1 > Q_2$; if $D_{Q_1} < D_{Q_2}$, then $Q_1 < Q_2$; if $D_{Q_1} = D_{Q_2}$, then $Q_1 = Q_2$.

**Definition 5[12]:** Let $Q_1 = <h_1, g_1> = <\{u_1^{(i)}(x)\}, \{v_1^{(i)}(x)\} >$ and $Q_2 = <h_2, g_2> < \{u_2^{(i)}(x)\}, \{v_2^{(i)}(x)\} >$ are two q-ROHFNs, then the *Hamming* distance between $h_1$ and $h_2$ as shown in formula (4).

$$D(Q_1, Q_2) = \frac{1}{2l}\sum_{i=1}^{l}\left\{\left|u_1^{(i)}(x) - u_2^{(i)}(x)\right| + \left|v_1^{(i)}(x) - v_2^{(i)}(x)\right| + \left|\pi_1^{(i)}(x) - \pi_2^{(i)}(x)\right|\right\}, \quad (4)$$

## 3 q-ROHFPR and its additive consistency

In this section, q-ROHFPR is defined and its properties are investigated, respectively. And additively consistent q-ROHFPR is proposed.

### 3.1 q-ROHFPR

When experts make decision on a set of alternatives $X = \{x_1, x_2, \ldots, x_n\}$, they use q-ROHFNs to express the evaluation value of the preference relation between two alternatives, so that the preference relation of the alternative can be expressed more flexibly. Inspired by q-ROHFS[12] and q-ROFPR[32], the definition of q-ROHFPR is proposed as shown in **Definition 6**.

**Definition 6:** A q-ROHFPR $A$ on $X$ is expressed by a matrix $A = (a_{ij})_{n\times n} \subseteq X \times X$, where $a_{ij} = <h_{ij}, g_{ij}> = <\{u_{ij}^{(s)}\}, \{v_{ij}^{(s)}\}>$ is a q-ROHFN for all $i, j = 1, 2, \ldots, n$, $h_{ij}$ characterizes the set of all possible membership degrees of $x_i \in X$, and $g_{ij}$ represents the set of all possible non-membership degrees of $x_i \in X$, with conditions:

$$\begin{cases} 0 \leq u_{ij}^{(s)} \leq 1 \\ 0 \leq v_{ij}^{(s)} \leq 1 \\ 0 \leq \left(u_{ij}^{(s)}\right)^q + \left(v_{ij}^{(s)}\right)^q \leq 1, \text{ for } q \geq 1, s = 1, 2, \ldots, n, i, j = 1, 2, \ldots, n \\ u_{ii} = v_{ii} = \{2^{-1/q}\} \\ u_{ij}^{(s)} = v_{ji}^{(s)}, u_{ji}^{(s)} = v_{ij}^{(s)} \end{cases} \quad (5)$$

In formula (5), $u_{ij}^{(s)} \in h_{ij}, v_{ij}^{(s)} \in g_{ij}$, $u_{ij}^{(s)}$ and $v_{ij}^{(s)}$ represent the value of the s-th element in $h_{ij}$ and $g_{ij}$, respectively. $u_{ij}^{(s)}$ represents the degree to which $x_i$ is preferred to $x_j$, $v_{ij}^{(s)}$ depicts the degree to which $x_j$ is preferred to $x_i$. If $u_{ij}^{(s)} = v_{ij}^{(s)} = 2^{-1/q}$, there is no difference between $x_i$ and $x_j$. The elements in $h_{ij}$ and $g_{ij}$ are arranged from small to large, and $l$ represents the number of elements in $h_{ij}$ and $g_{ij}$. The fuzzy preference relation of n alternatives can be expressed as a n*n matrix of q-ROHFPR, as shown in formula (6).

$$A = (a_{ij})_{n \times n} = \begin{pmatrix} <\{2^{-1/q}\},\{2^{-1/q}\}> & <\{u_{12}^{(1)},u_{12}^{(2)},\ldots,u_{12}^{(l)}\},\{v_{12}^{(1)},v_{12}^{(2)},\ldots,v_{12}^{(l)}\}> & \cdots & <\{u_{14}^{(1)},u_{14}^{(2)},\ldots,u_{14}^{(l)}\},\{v_{14}^{(1)},v_{14}^{(2)},\ldots,v_{14}^{(l)}\}> \\ <\{v_{12}^{(1)},v_{12}^{(2)},\ldots,v_{12}^{(l)}\},\{u_{12}^{(1)},u_{12}^{(2)},\ldots,u_{12}^{(l)}\}> & <\{2^{-1/q}\},\{2^{-1/q}\}> & \cdots & <\{u_{24}^{(1)},u_{24}^{(2)},\ldots,u_{24}^{(l)}\},\{v_{24}^{(1)},v_{24}^{(2)},\ldots,v_{24}^{(l)}\}> \\ \vdots & \vdots & \vdots & <\{u_{34}^{(1)},u_{34}^{(2)},\ldots,u_{34}^{(l)}\},\{v_{34}^{(1)},v_{34}^{(2)},\ldots,v_{34}^{(l)}\}> \\ <\{v_{14}^{(1)},v_{14}^{(2)},\ldots,v_{14}^{(l)}\},\{u_{14}^{(1)},u_{14}^{(2)},\ldots,u_{14}^{(l)}\}> & <\{v_{24}^{(1)},v_{24}^{(2)},\ldots,v_{24}^{(l)}\},\{u_{24}^{(1)},u_{24}^{(2)},\ldots,u_{24}^{(l)}\}> & \cdots & <\{2^{-1/q}\},\{2^{-1/q}\}> \end{pmatrix}$$
(6)

The q-ROHFPR matrix not only expresses the direction of preference among several alternatives, but also expresses the strength of preference. In order to ensure the consistency of experts' opinions, this definition limits the number of preference membership degrees of each hesitant fuzzy number in $A$ to be equal to the number of preference non-membership degree elements, that is, in formula (6), each q-ROHFNs contains $l$ preference for membership and non-membership, because when there are and only two alternatives, i.e. $x_i$ and $x_j$ are compared, $l$ is not required, but when there are three or more alternatives to compare, if it means $x_i, x_j, x_k$, When the preference membership degree and non-membership degree of q-ROHFNs of them are not equal, it is easy to deduce: $u_{ij}^{(s)} + u_{jk}^{(s)} + u_{ki}^{(s)}$ and $u_{ik}^{(s)} + u_{kj}^{(s)} + u_{ji}^{(s)}$ cannot be kept equal. According to the value of matrix $A$, it is easy to judge the preference relation between the two alternatives, but it is impossible to directly judge the relation between multiple alternatives. When each q-ROHFNs in the preference relation matrix $A$ contains $l$ preference membership degrees and non-membership degrees are equal, inspired by [32], the following 6 properties can be used to determine the relationship between multiple elements in the q-ROHFPR matrix $A$.

(1) Triangle condition. The q-ROHFPR $A$ satisfies the triangle condition if $a_{ik} \oplus a_{kj} \geq a_{ij}$ holds for $\forall a_{ij}, a_{ik}, a_{kj} \in \Omega, i, j, k = 1, 2, \ldots, n, q > 1$. Let alternatives $x_i, x_j, x_k$ be three vertices of a triangle with length sides $a_{ij}, a_{ik}$ and $a_{kj}$, then the length of vertices $x_i, x_j$ should be smaller than the summation of the lengths with respect to the vertices $x_i, x_k$ and $x_k, x_j$.

(2) Weak transitivity property. The q-ROHFPR $A$ satisfies the weak transitivity property if $a_{ik} \geq < \{2^{-\frac{1}{q}}\}, \{2^{-\frac{1}{q}}\} >$, $a_{kj} \geq < \{2^{-\frac{1}{q}}\}, \{2^{-\frac{1}{q}}\} > \Longrightarrow a_{ij} \geq < \{2^{-\frac{1}{q}}\}, \{2^{-\frac{1}{q}}\} >$ holds for $\forall a_{ij}, a_{ik}, a_{kj} \in \Omega, i, j, k = 1, 2, \ldots, n, q > 1$. This property shows that if alternatives $x_i$ is preferred to $x_k$, and $x_k$ is preferred to $x_j$, then alternative $x_i$ should be preferred to alternative $x_j$. It is noted that the weak transitivity property is the minimum requirement that a consistent q-ROHFPR should verify.

(3) Max-min transitivity property. The q-ROHFPR $A$ satisfies the max-min transitivity property if $a_{ij} \geq \min\{a_{ik}, a_{kj}\}$ holds for $\forall a_{ij}, a_{ik}, a_{kj} \in \Omega, i, j, k = 1,2, \ldots, n, q > 1$. It means that the q-ROHFN obtained from a direct comparison between two alternatives should be greater than or equal to the minimum partial value obtained by comparing alternatives with an intermediate one.

(4) Max-max transitivity property. The q-ROHFPR $A$ satisfies the max-max transitivity property if $a_{ij} \geq \max\{a_{ik}, a_{kj}\}$ holds for $\forall a_{ij}, a_{ik}, a_{kj} \in \Omega, i, j, k = 1,2, \ldots, n, q > 1$. It denotes that the q-ROHFN obtained from a direct comparison between two alternatives should be greater than or equal to the maximum partial values obtained by comparing alternatives with an intermediate one.

(5) Restricted max-min transitivity property. The q-ROHFPR $A$ satisfies the restricted max-min transitivity property if $a_{ik} \geq <\{2^{-\frac{1}{q}}\}, \{2^{-\frac{1}{q}}\}>, a_{kj} \geq <\{2^{-\frac{1}{q}}\}, \{2^{-\frac{1}{q}}\}> \Longrightarrow a_{ij} \geq \min\{h_{ik}, h_{kj}\}$ holds for $\forall a_{ij}, a_{ik}, a_{kj} \in \Omega, i, j, k = 1,2, \ldots, n, q > 1$. It presents that if alternative $x_i$ is preferred to $x_k$ with a q-ROHFN $a_{ik}$, and $x_k$ is preferred to $x_j$ with a q-ROHFN $a_{kj}$, then $x_i$ should be preferred to $x_j$ with a q-ROHFN $a_{ij}$ which is greater than or equal to the minimum value of $a_{ik}$ and $a_{kj}$. The equality holds when there is indifference between at least two of the three alternatives.

(6) Restricted max-max transitivity property. The q-ROHFPR $A$ satisfies the restricted max-max transitivity property if $a_{ik} \geq <\{2^{-\frac{1}{q}}\}, \{2^{-\frac{1}{q}}\}>, a_{kj} \geq <\{2^{-\frac{1}{q}}\}, \{2^{-\frac{1}{q}}\}> \Longrightarrow a_{ij} \geq \max\{a_{ik}, a_{kj}\}$ holds for $\forall a_{ij}, a_{ik}, a_{kj} \in \Omega, i, j, k = 1,2, \ldots, n, q > 1$. It expresses that if alternative $x_i$ is preferred to $x_k$ with a q-ROHFN $a_{ik}$, and $x_k$ is preferred to $x_j$ with a q-ROHFN $a_{kj}$, then $x_i$ should be preferred to $x_j$ with a q-ROHFN $a_{ij}$ which is greater than or equal to the minimum value of $a_{ik}$ and $a_{kj}$. Similarly, the equality holds when there is indifference between at least two of the three alternatives.

Regarding the above six properties, the triangular inequality property in the first (1) is the basic relationship between the three elements, and the weak transitivity in the (2) can be the basic condition for meeting the consistency.

### 3.2 Additively consistent q-ROHFPR

It is necessary to judge and verify that the preference relationship given by DMs is consistent to reflect the logic of their rationality of the judgment. Lack of consistency in decision-making may lead to inconsistent conclusions. In the decision-making process, the consistency condition is usually to meet the transitivity, that is, there are any three alternatives: $x_i, x_j, x_k$, if $x_i$ is better than $x_j$, $x_j$ is better than $x_k$, then $x_i$ should be better than $x_k$. According to the judgment value given by DMs, there are usually additive transitivity and multiplicative transitivity. Inspired by Zhang and Liao[32], a definition of additive consistency of q-ROHFPR based on additive transitivity was proposed.

**Definition 7:** A q-ROHFPR $A = (a_{ij})_{n \times n}$ is called additively consistent for all $i, j, k \in N, s \in \{1, 2, \ldots, l\}$, if it satisfies the following additive transitivity:

$$\begin{cases} (u_{ij}^{(s)})^q + (u_{jk}^{(s)})^q + (u_{ki}^{(s)})^q = (u_{ik}^{(s)})^q + (u_{kj}^{(s)})^q + (u_{ji}^{(s)})^q \\ (v_{ij}^{(s)})^q + (v_{jk}^{(s)})^q + (v_{ki}^{(s)})^q = (v_{ik}^{(s)})^q + (v_{kj}^{(s)})^q + (v_{ji}^{(s)})^q \end{cases}, \quad (7)$$

In formula (11), $a_{ij} = <\{u_{ij}^{(s)}\}, \{v_{ij}^{(s)}\}>$, $0 \leq (u_{ij}^{(s)})^q + (v_{ij}^{(s)})^q \leq 1, q \geq 1$, $u_{ij}^{(s)}$ and $v_{ij}^{(s)}$ are respectively the $s$th membership element and the $s$th non-membership element in the hesitant fuzzy element (HFE) $a_{ij}$.

This definition is based on FPR $B = (b_{ij})_{n \times n}$ additive transitivity formula $b_{ij} + b_{jk} + b_{ki} = b_{ik} + b_{kj} + b_{ji}$ proposed by Tanino[33], for all $i, j, k \in N$. In addition, when the number of elements in $a_{ij}$ is 1, that is, when $\{u_{ij}^{(s)}\}$ and $\{v_{ij}^{(s)}\}$ in $<\{u_{ij}^{(s)}\}, \{v_{ij}^{(s)}\}>$ degenerate into a real number $u_{ij}$ and $u_{ij}$, respectively, it can be known from formula (11) that for all $i, j, k \in N$, holds $\begin{cases} u_{ij} + u_{jk} + u_{ki} = u_{ik} + u_{kj} + u_{ji} \\ v_{ij} + v_{jk} + v_{ki} = v_{ik} + v_{kj} + v_{ji} \end{cases}$, which is equal to the additive consistency condition of q-ROFPR in [32].

Based on **Definition 6** and **Definition 7**, a new additive consistency definition of q-ROHFPRs is proposed, as shown in **Definition 8**.

**Definition 8:** A q-ROHFPR $A = (a_{ij})_{n \times n}$ is called additively consistent with $a_{ij} = <\{u_{ij}^{(s)}\}, \{v_{ij}^{(s)}\}>$, $0 \leq (u_{ij}^{(s)})^q + (v_{ij}^{(s)})^q \leq 1$ for all $i, j = 1, 2, \ldots, n, s \in \{1, 2, \ldots, l\}, q \geq 1$, if it satisfies the following formula (8).

$$(u_{ij}^{(s)})^q + (u_{jk}^{(s)})^q - (u_{ik}^{(s)})^q = (v_{ij}^{(s)})^q + (v_{jk}^{(s)})^q - (v_{ik}^{(s)})^q, i < j < k, s = 1, 2, \ldots, l, \quad (8)$$

According to the **definition 6**, $u_{ij}^{(s)} = v_{ji}^{(s)}$, formula (11) can be reduced to $(u_{ij}^{(s)})^q + (u_{jk}^{(s)})^q + (u_{ki}^{(s)})^q = (v_{ij}^{(s)})^q + (v_{jk}^{(s)})^q + (v_{ki}^{(s)})^q$ for all $i<j<k$, obviously $u_{ij}^{(s)}, u_{jk}^{(s)}$ are located in the upper triangle of $A$, and $u_{ki}^{(s)}$, $u_{ji}^{(s)}$ are located in the lower triangle of $A$, for the convenience of calculation, replace $v_{ik}^{(s)}$ with $u_{ki}^{(s)}$, and similarly replace $u_{ik}^{(s)}$ with $v_{ki}^{(s)}$ to obtain formula (8). Due to the high complexity of calculating transfer consistency in the solution process, it is not suitable for adjustment. Therefore, on the basis of formula (8), the definition of additive consistency index is proposed.

**Definition 9:** The consistency index CI(A) of a q-ROHFPR $A = (a_{ij})_{n\times n} = (<\{u_{ij}^{(s)}\},\{v_{ij}^{(s)}\}>)_{n\times n}$ is defined as

$$CI(A)=\frac{1}{l}\frac{2}{n(n-1)(n-2)}\sum_{1\leq i<j<k}^{n}\sum_{s=1}^{l}\left|(u_{ij}^{(s)})^q+(u_{jk}^{(s)})^q+(v_{ik}^{(s)})^q-\left(v_{ij}^{(s)}\right)^q-(v_{jk}^{(s)})^q-(u_{ik}^{(s)})^q\right|, \quad (9)$$

where $0 \leq \left(u_{ij}^{(s)}\right)^q + \left(v_{ij}^{(s)}\right)^q \leq 1$, $i,j = 1,2,\ldots,n$, $s \in \{1,2,\ldots,l\}$, $q \geq 1$. In formula (13), since $0 \leq \left(u_{ij}^{(s)}\right)^q + \left(v_{ij}^{(s)}\right)^q \leq 1$, it is easy to deduce CI(A) $\in [0,1]$. At the same time, when the value of the consistency index CI is smaller, the consistency level of the q-ROHFPR $A$ is higher, when CI(A) = 0, q-ROHFPR A is completely consistent. However, in actual decision-making problems, the cost of solving complete consistency q-ROHFPR is too high. For this reason, the consistency index is often used to determine the acceptable consistency of q-ROHFPR, that is, when CI is less than or equal to a threshold $\overline{CI}$, consider q-ROHFPR is acceptably consistent. For this reason, the acceptably additive consistent q-ROHFPR is defined, as shown in **Definition 10**.

**Definition 10:** The q-ROHFPR $A$ is called to be acceptably additive consistent if CI(A) $\leq \overline{CI}$; otherwise, $A$ is called unacceptably additive consistent. When CI(A) = 0, $A$ is completely consistent, where $\overline{CI}$ is the consistency threshold and $\overline{CI}(A) \in [0,1]$.

## 4 Consistency models for q-ROHFPRs

Several goal programming models are developed to repair individual consistency and derived q-rung orthopair hesitant fuzzy weights from the q-ROHFPRs.

### 4.1 Acceptably consistent goal programming model of q-ROHFPRs

In the actual DM process, the degree of consistency of a q-ROHFPR represents its rationality, but due to the uncertainty of decision-making information, it is difficult for DMs to directly provide completely additive consistency or acceptable additive consistency q-ROHFPRs. If the consistency of q-ROHFPRs are poor, the decision-making results derived from the q-ROHFPR may be suspicious and unsatisfactory. In order to solve this situation, several goal programming models have been developed to improve the consistency of given q-ROHFPRs, and an acceptable additive consistency q-ROHFPR $A = (a_{ij})_{n \times n}$ is derived from the unacceptable additive consistency q-ROHFPR. Let $A$ be the unacceptable additive consistency q-ROHFPR, and $\tilde{A}$ be the adjusted q-ROHFPR of $A$. The Manhattan distance of $A = (a_{ij})_{n \times n}$ and $\tilde{A} = (\tilde{a}_{ij})_{n \times n}$ are defined.

**Definition 11:** Let $A = (a_{ij})_{n \times n}$ and $\tilde{A} = (\tilde{a}_{ij})_{n \times n}$ be two q-ROHFPRs. For any $i, j \in N, s = 1, 2, \ldots, l$, the Manhattan distance between $A$ and $\tilde{A}$ is shown in formula (10).

$$D(A, \tilde{A}) = \frac{1}{2l} \frac{1}{(n-1)(n-2)} \sum_{1 \leq i < j \leq n}^{n} \sum_{s=1}^{l} \left\{ \left| u_{ij}^{(s)} - \tilde{u}_{ij}^{(s)} \right| + \left| v_{ij}^{(s)} - \tilde{v}_{ij}^{(s)} \right| \right\}, \tag{10}$$

In formula (10), $a_{ij} = <\{u_{ij}^{(s)}\}, \{v_{ij}^{(s)}\}>$, $\tilde{a}_{ij} = <\{\tilde{u}_{ij}^{(s)}\}, \{\tilde{v}_{ij}^{(s)}\}>$. Because of $0 \leq |u_{ij}^{(s)} - \tilde{u}_{ij}^{(s)}|, |v_{ij}^{(s)} - \tilde{v}_{ij}^{(s)}| \leq 1$, we set the coefficient to $\frac{1}{2l} \frac{1}{(n-1)(n-2)}$ to limit the restriction of $D(A, \tilde{A})$ in $[0,1]$. According to equation (10), $D(A, \tilde{A}) = D(\tilde{A}, A)$. In addition, it can be seen that the smaller the distance measure $D(A, \tilde{A})$ is, the closer the q-ROHFPRs are.

The model proposed in this paper aims to minimize the overall adjustment of q-ROHFPR. In this case, all elements in q-ROHFPR should have the same importance. Therefore, the deviation between two q-ROHFPRs $A$ and $\tilde{A}$ is proposed based on Definition 12.

**Definition 12:** Suppose $A = (a_{ij})_{n \times n}$ and $\tilde{A} = (\tilde{a}_{ij})_{n \times n}$ are two q-ROHFPRs. For any $i, j \in N, s = 1, 2, \ldots, l$, the deviation between $A$ and $\tilde{A}$ is shown in formula (11).

$$dev(A, \tilde{A}) = \sum_{1 \leq i < j \leq n}^{n} \sum_{s=1}^{l} \left\{ \left| u_{ij}^{(s)} - \tilde{u}_{ij}^{(s)} \right| + \left| v_{ij}^{(s)} - \tilde{v}_{ij}^{(s)} \right| \right\}, \tag{11}$$

In formula (11), $a_{ij} = <\{u_{ij}^{(s)}\}, \{v_{ij}^{(s)}\}>$, $\tilde{a}_{ij} = <\{\tilde{u}_{ij}^{(s)}\}, \{\tilde{v}_{ij}^{(s)}\}>$.

For the unacceptable consistency q-ROHFPR $A$, it is necessary to obtain a new q-ROHFPR that satisfies the acceptable additive consistency and retains the original preference information as much as possible. From the definition of $\tilde{u}_{ij}^{(s)} = \tilde{v}_{ji}^{(s)}$ and $\tilde{u}_{ji}^{(s)} = \tilde{v}_{ij}^{(s)}$ in the generalized hesitation fuzzy preference relationship, it can be known that the distance of the upper triangle element of Matrix A is equal to the distance of the lower triangle element, so we only need to request the adjustment element of the upper triangle element of matrix A. From the perspective of goal planning, in order to achieve this goal, given a q-ROHFPR $A$ and a predefined consistency threshold $\overline{CI}$, under the constraints of acceptable additive consistency, by minimizing the new q- The absolute deviation between q-ROHFPR and original q-ROHFPR $A$, the mathematical programming model is established as shown in $(M-1)$.

$$\min dev(A, \tilde{A})$$
$$s.t. \begin{cases} CI(\tilde{A}) \leq \overline{CI} \\ \tilde{A} \quad is \quad a \quad q-ROHFPR. \end{cases} \quad \textbf{(M-1)}$$

In $(M-1)$, $A$ is the q-ROHFPR matrix given by the expert, and $\overline{CI}(0 \leq \overline{CI} \leq 1)$ is the predefined consistency threshold, and the values of both are fixed. $(M-1)$ is a nonlinear programming model. In the model $(M-1)$, $\tilde{A}$ is a q-ROHFPR, $CI(\tilde{A}) \leq \overline{CI}$ are two constraints, and $dev(A, \tilde{A})$ is the objective function. The first constraint in the model $(M-1)$ guarantees that the matrix $\tilde{A}$ obtained by the model $(M-1)$ wreaks the requirements of acceptable additive consistency, while the second constraint is to ensure that $\tilde{A}$ is an optimized q-ROHFPR matrix. The model $(M-1)$ aims to determine a new q-ROHFPR $\tilde{A}$ with the smallest distance from the original q-ROHFPR $A$ under the premise of satisfying these two constraints. Based on the above analysis, we can guarantee that the absolute error obtained through the model $(M-1)$ will be the smallest, and the consistency index of the new q-ROHFPR $\tilde{A}$ obtained through the model $(M-1)$ satisfies the prescribed threshold. Therefore, The $\tilde{A}$ obtained by the model $(M-1)$ is q-ROHFPR with acceptable consistency. Substituting (13) and (15) into $(M-1)$, $(M-1)$ can be expressed as $(M-2)$.

$$\min \sum_{i=1}^{n-1} \sum_{j=i+1}^{n} \sum_{s=1}^{l} \left( \left| u_{ij}^{(s)} - \tilde{u}_{ij}^{(s)} \right| + \left| v_{ij}^{(s)} - \tilde{v}_{ij}^{(s)} \right| \right)$$

$$s.t. \begin{cases} \dfrac{1}{l} \dfrac{2}{n(n-1)(n-2)} \sum_{1 \le i < j < n} \sum_{s=1}^{l} \{ \left| \tilde{u}_{ij}^{(s)} + \tilde{u}_{jk}^{(s)} + \tilde{v}_{ik}^{(s)} - \tilde{v}_{ij}^{(s)} - \tilde{v}_{jk}^{(s)} - \tilde{u}_{ik}^{(s)} \right| \} \le \overline{CI}, \\ \tilde{u}_{ii}^{(s)} = \tilde{v}_{ii}^{(s)} = (\tfrac{1}{2})^q, \tilde{u}_{ij}^{(s)} = \tilde{v}_{ji}^{(s)}, \\ \tilde{u}_{ij}^1 < \tilde{u}_{ij}^2 < \ldots < \tilde{u}_{ij}^s, \tilde{v}_{ij}^1 < \tilde{v}_{ij}^2 < \ldots < \tilde{v}_{ij}^s, \\ 0 \le \tilde{u}_{ij}^s, \tilde{v}_{ij}^s \le 1, \\ 0 \le (\tilde{u}_{ij}^s)^q + (\tilde{v}_{ij}^s)^q \le 1, \\ i, j, k = 1, 2, \ldots, n;\ s = 1, 2, \ldots, l \end{cases} \quad \text{(M-2)}$$

By solving the model $(M-2)$, the optimal solution $\tilde{q}_{ij} = <\{\tilde{u}_{ij}^{(s)}\}, \{\tilde{v}_{ij}^{(s)}\}>$ is the acceptable additive consistency. The upper triangular matrix element of q-ROHFPR $\tilde{A}$, through the definition 9, can be completed to get the complete q-ROHFPR $\tilde{A}$, which is the goal we are seeking, and we can use it to complete the subsequent decision-making.

### 4.2 Priority model of alternatives

In this subsection, we will construct a priority vectors of alternatives based on model $(M-2)$. Let $w = (w_1, w_2, \ldots, w_n)^T$ be the priority vector of q-ROHFPR $A = (a_{ij})_{n \times n}$, where $a_{ij} = <\{u_{ij}^{(s)}\}, \{v_{ij}^{(s)}\}>$, $w_i = <\{w_{iu}^s\}, \{w_{iv}^s\}>$, $0 \le w_{iu}^s, w_{iv}^s \le 1, 0 \le (w_{iu}^s)^q + (w_{iv}^s)^q \le 1$. If $q-ROHF$ priority vector $w = (w_1, w_2, \ldots, w_n)^T$ satisfies $\sum_{j=1, j \ne i}^{n} (w_{iu}^s)^q \le (w_{iv}^s)^q$ and $(w_{iu}^s)^q + n - 2 \ge \sum_{j=1, j \ne i}^{n} (w_{iv}^s)^q$, then w is a normalized vector. The q-ROHFPR $\tilde{A} = (\tilde{a}_{ij})_{n \times n}$ is constructed to solve the q-ROHF priority based on adding consistency q-ROHFPR proposed by Zhang[34].

**Definition 13：** Let $A = (a_{ij})_{n \times n}$ is a q-ROHFPR, where $a_{ij} = <\{u_{ij}^s\}, \{v_{ij}^s\}>$, if there is a normalized q-ROHF priority weight vector $w = (w_1, w_2, \ldots, w_n)^T$ and $w_i = <\{w_{iu}^s\}, \{w_{iv}^s\}>$ satisfies formula (12), then $A$ is q-ROHFPR with consistency.

$$a_{ij} = <\{u_{ij}^s\}, \{v_{ij}^s\}> = \begin{cases} \langle \{\sqrt[q]{0.5}\}, \{\sqrt[q]{0.5}\} \rangle & \text{if } i = j \\ \langle \{\sqrt[q]{0.5(w_{iu}^s)^q + 0.5(w_{jv}^s)^q}\}, \{\sqrt[q]{0.5(w_{iv}^s)^q + 0.5(w_{ju}^s)^q}\} \rangle & \text{if } i \ne j \end{cases}, \quad (12)$$

where $q \ge 1, i, j = 1, 2, \ldots, n, s = 1, 2, \ldots, l$.

According to formula (12), if q-ROHFPR $A = (a_{ij})_{n \times n}$ is plus consistency, where $a_{ij} = <\{u_{ij}^s\}, \{v_{ij}^s\}>, i < j = 1,2,\ldots,n, s = 1,2,\ldots,l$, we can get:

$$\begin{cases} u_{ij}^s = \sqrt[q]{0.5(w_{iu}^s)^q + 0.5(w_{jv}^s)^q} \\ v_{ij}^s = \sqrt[q]{0.5(w_{iv}^s)^q + 0.5(w_{ju}^s)^q} \end{cases} \quad (13)$$

Next, To solve the priority weight of q-ROHFPR, two variables $\lambda$ and $\theta$ are introduced to measure the deviation between the degree of membership, the degree of non-membership and the weight respectively:

$$\begin{cases} \lambda_{ij}^s = u_{ij}^s - \sqrt[q]{0.5(w_{iu}^s)^q + 0.5(w_{jv}^s)^q} \\ \theta_{ij}^s = v_{ij}^s - \sqrt[q]{0.5(w_{iv}^s)^q + 0.5(w_{ju}^s)^q} \end{cases}, \text{for } i,j = 1,2,\ldots,n, s = 1,2,\ldots,l \quad (14)$$

According to formula (14), if $\lambda_{ij}^s + \theta_{ij}^s = 0$, then q-ROHFPR satisfies additive consistency. The smaller the sum of $|\lambda_{ij}^s|$ and $|\theta_{ij}^s|$, the higher the consistency level of q-ROHFPR. Therefore, the model $(M-3)$ is constructed by minimizing the value of $(|\lambda_{ij}^s| + |\theta_{ij}^s|)$:

$$\min \sum_{i=1}^{n-1} \sum_{j=i+1}^{n} \sum_{s=1}^{l} |\lambda_{ij}^s| + |\theta_{ij}^s|$$

$$\text{s.t.} \begin{cases} \lambda_{ij}^s = u_{ij}^s - \sqrt[q]{0.5(w_{iu}^s)^q + 0.5(w_{jv}^s)^q}, \\ \theta_{ij}^s = v_{ij}^s - \sqrt[q]{0.5(w_{iv}^s)^q + 0.5(w_{ju}^s)^q}, \\ 0 \leq w_{iu}^s, w_{iv}^s \leq 1, \\ 0 \leq (w_{iu}^s)^q + (w_{iv}^s)^q \leq 1, \\ \sum_{j=1, j \neq i}^{n} (w_{ju}^s)^q \leq (w_{iv}^s)^q, \\ (w_{iu}^s)^q + n - 2 \geq \sum_{j=1, j \neq i}^{n} (w_{jv}^s)^q, \\ i = 1,2,\ldots,n; s = 1,2,\ldots,l \end{cases} \quad \text{(M-3)}$$

Let $|\lambda_{ij}^s| = \lambda_{ij}^{s+} + \lambda_{ij}^{s-}$, $|\theta_{ij}^s| = \theta_{ij}^{s+} + \theta_{ij}^{s-}$, satisfy $\lambda_{ij}^{s+} \times \lambda_{ij}^{s-} = 0$, $\theta_{ij}^{s+} \times \theta_{ij}^{s-} = 0$, $i = 1,2,\ldots,n-1$, $j = i+1$, $s = 1,2,\ldots,l$, the above model can be rewritten as: $(M-4)$

$$\min \sum_{i=1}^{n-1} \sum_{j=i+1}^{n} \sum_{s=1}^{l} \lambda_{ij}^{s+} + \lambda_{ij}^{s-} + \theta_{ij}^{s+} + \theta_{ij}^{s-}$$

$$s.t. \begin{cases} u_{ij}^s - \sqrt[q]{0.5(w_{iu}^s)^q + 0.5(w_{jv}^s)^q} - \lambda_{ij}^{s+} + \lambda_{ij}^{s-} = 0, \\ v_{ij}^s - \sqrt[q]{0.5(w_{iv}^s)^q + 0.5(w_{ju}^s)^q} - \theta_{ij}^{s+} + \theta_{ij}^{s-} = 0, \\ 0 \leq w_{iu}^s, w_{iv}^s \leq 1, \\ 0 \leq (w_{iu}^s)^q + (w_{iv}^s)^q \leq 1, \\ \sum_{j=1, j \neq i}^{n} (w_{ju}^s)^q \leq (w_{iv}^s)^q, \\ (w_{iu}^s)^q + n - 2 \geq \sum_{j=1, j \neq i}^{n} (w_{jv}^s)^q, \\ \lambda_{ij}^{s+}, \lambda_{ij}^{s-}, \theta_{ij}^{s+}, \theta_{ij}^{s-} \geq 0, \lambda_{ij}^{s+} \times \lambda_{ij}^{s-} = 0, \theta_{ij}^{s+} \times \theta_{ij}^{s-} = 0 \\ i = 1, 2, \ldots, n; s = 1, 2, \ldots, l \end{cases} \quad (M-4)$$

By solving the model $(M-4)$, the q-ROHF priority weight $w = (w_1, w_2, \ldots, w_n)^T$ can be obtained.

## 5 Group decision making

A consensus index is proposed, a new method is proposed to solve the GDM with q-ROHFPRs.

### 5.1 A consensus index for q-ROHFPRs

Let $X = \{x_1, x_2, \ldots, x_n\}$ be a set of alternatives and $D = \{d_1, d_2, \ldots, d_m\}$ be a set of DMs whose weight vector is $\lambda = (\lambda_1, \lambda_2, \ldots, \lambda_m)^T$, where $\lambda_t$ is the weight of DM $d_t$ satisfying $\lambda_t \geq 0, \sum_{t=1}^{m} \lambda_t = 1, (t = 1, 2, \ldots, m)$. Let q-ROHFPRs $A_t = (a_{ijt})_{n \times n}$ $(t = 1, 2, \ldots, m)$ with $a_{ijt} = (\mu_{ijt}, v_{ijt})$ be $m$ individual decision matrices, where $\mu_{ijt} = \{\mu_{ijt}^1, \mu_{ijt}^2, \ldots, \mu_{ijt}^l\}$, $v_{ijt} = \{v_{ijt}^1, v_{ijt}^2, \ldots, v_{ijt}^l\}$ and $a_{ijt}$ is an attribute value, given by the expert $d_t$. If some individual q-ROHFPRs are unacceptable consistent, it can be repaired by (M-2). Consensus refers to the unanimity of individual opinions that reflects the option of the group. To obtain the representative ranking, it is necessary to give common understanding analysis on group minds. In GDM problem, when individual q-ROHFPRs are additively consistent, it is also expected that aggregated q-ROHFPRs are additively consistent. Therefore, inspired by the AWA operator in Ref.[30], an aggregation operator is proposed, as shown in **Definition 14**.

**Definition 14:** : Let $A_t = (a_{ijt})_{n \times n}$ be the q-ROHFPRs of acceptable additive consistency on the alternatives set $X = \{x_1, x_2, \ldots, x_n\}$, where $a_{ijt} = (\{\mu_{ijt}^1, \mu_{ijt}^2, \ldots, \mu_{ijt}^l\}, \{v_{ijt}^1, v_{ijt}^2, \ldots, v_{ijt}^l\})$, $t = 1, 2, \ldots, m$. $A = (a_{ij})_{n \times n}$ is the aggregation matrix of $A_t$, where $a_{ij} = (\{\mu_{ij}^1, \mu_{ij}^2, \ldots, \mu_{ij}^l\}, \{v_{ij}^1, v_{ij}^2, \ldots, v_{ij}^l\})$, which can be aggregated by the following operator, as shown in formula (15).

$$a_{ij} = \sum_{t=1}^{m} \lambda_t a_{ijt}, \; t = 1, 2, \ldots, m, \; i, j = 1, 2, \ldots, n \tag{15}$$

namely

$$a_{ij} = (\mu_{ij}, v_{ij}) = (\{\mu_{ij}^1, \mu_{ij}^2, \ldots, \mu_{ij}^l\}, \{v_{ij}^1, v_{ij}^2, \ldots, v_{ij}^l\})$$
$$= (\{(\sum_{k=1}^{m} \lambda_t (\mu_{ijt}^1)^q)^{1/q}, (\sum_{t=1}^{m} \lambda_t (\mu_{ijt}^2)^q)^{1/q}, \ldots,$$
$$(\sum_{t=1}^{m} \lambda_t (\mu_{ijt}^l)^q)^{1/q}\}, \{(\sum_{t=1}^{m} \lambda_t (v_{ijt}^1)^q)^{1/q}, (\sum_{t=1}^{m} \lambda_t (v_{ijt}^2)^q)^{1/q}, \ldots, (\sum_{t=1}^{m} \lambda_t (v_{ijt}^l)^q)^{1/q}\})$$

(16)

In formula (15) and (16), $t = 1, 2, \ldots, m$, $i, j = 1, 2, \ldots, n$, $\lambda = (\lambda_1, \lambda_2, \ldots, \lambda_m)^T$ be the weight vector of DMs, where $0 \leq \lambda_t \leq 1$, $\sum_{t=1}^{m} \lambda_t = 1$. $a_{ijt} = (\{\mu_{ijt}^1, \mu_{ijt}^2, \ldots, \mu_{ijt}^l\}, \{v_{ijt}^1, v_{ijt}^2, \ldots, v_{ijt}^l\})$, $0 \leq (\mu_{ijt}^l)^q + (v_{ijt}^l)^q \leq 1$, $q \geq 1$. $A$ is also a q-ROHFPR.

**Theorem 1**: If all individual q-ROHFPRs satisfy the acceptable additive consistency, then the aggregated q-ROHFPR also satisfies the acceptable additive consistency.

**Proof.** Assume that all individual q-ROHFPRs $A_t = (a_{ijt})_{n \times n}$ $t = 1, 2, \ldots, m$ are of acceptably additive consistency, according to **Definition 10**, one can get $CI(A_t) \leq \overline{CI}$. To prove that the aggregated q-ROHFPR $A$ is acceptably additive consistent, we need to demonstrate that $CI(A) \leq \overline{CI}$.

From **Definition 9**, we have

$$CI(A) = \frac{1}{l} \frac{2}{n(n-1)(n-2)} \sum_{1 \leq i < j < k} \sum_{s=1}^{l} \left| (\mu_{ij}^s)^q + (\mu_{jk}^s)^q + (v_{ik}^s)^q - (v_{ij}^s)^q - (v_{jk}^s)^q - (\mu_{ik}^s)^q \right|$$

Let $\Phi = (\mu_{ij}^s)^q + (\mu_{jk}^s)^q + (v_{ik}^s)^q - (v_{ij}^s)^q - (v_{jk}^s)^q - (\mu_{ik}^s)^q$, combined with Eq.(16), it can be deduced

$$\Phi = \left(\mu_{ij}^{s}\right)^{q} + \left(\mu_{jk}^{s}\right)^{q} + \left(v_{ik}^{s}\right)^{q} - \left(v_{ij}^{s}\right)^{q} - \left(v_{jk}^{s}\right)^{q} - \left(\mu_{ik}^{s}\right)^{q}$$

$$= \left[\left(\sum_{t=1}^{m}\lambda_{t}\left(\mu_{ijt}^{s}\right)^{q}\right)^{1/q}\right]^{q} + \left[\left(\sum_{t=1}^{m}\lambda_{t}\left(\mu_{jkt}^{s}\right)^{q}\right)^{1/q}\right]^{q} + \left[\left(\sum_{t=1}^{m}\lambda_{t}\left(v_{ikt}^{s}\right)^{q}\right)^{1/q}\right]^{q} - \left[\left(\sum_{t=1}^{m}\lambda_{t}\left(v_{ijt}^{s}\right)^{q}\right)^{1/q}\right]^{q} - \left[\left(\sum_{t=1}^{m}\lambda_{t}\left(v_{jkt}^{s}\right)^{q}\right)^{1/q}\right]^{q} - \left[\left(\sum_{t=1}^{m}\lambda_{t}\left(\mu_{ikt}^{s}\right)^{q}\right)^{1/q}\right]^{q}$$

$$= \sum_{t=1}^{m}\lambda_{t}\left(\mu_{ijt}^{s}\right)^{q} + \sum_{t=1}^{m}\lambda_{t}\left(\mu_{jkt}^{s}\right)^{q} + \sum_{t=1}^{m}\lambda_{t}\left(v_{ikt}^{s}\right)^{q} - \sum_{t=1}^{m}\lambda_{t}\left(v_{ijt}^{s}\right)^{q} - \sum_{t=1}^{m}\lambda_{t}\left(v_{jkt}^{s}\right)^{q} - \sum_{t=1}^{m}\lambda_{t}\left(\mu_{ikt}^{s}\right)^{q}$$

$$= \lambda_{t}\left(\sum_{t=1}^{m}\left(\mu_{ijt}^{s}\right)^{q} + \sum_{t=1}^{m}\left(\mu_{jkt}^{s}\right)^{q} + \sum_{t=1}^{m}\left(v_{ikt}^{s}\right)^{q} - \sum_{t=1}^{m}\left(v_{ijt}^{s}\right)^{q} - \sum_{t=1}^{m}\left(v_{jkt}^{s}\right)^{q} - \sum_{t=1}^{m}\left(\mu_{ikt}^{s}\right)^{q}\right)$$

$$= \lambda_{t}\sum_{t=1}^{m}\left(\left(\mu_{ijt}^{s}\right)^{q} + \left(\mu_{jkt}^{s}\right)^{q} + \left(v_{ikt}^{s}\right)^{q} - \left(v_{ijt}^{s}\right)^{q} - \left(v_{jkt}^{s}\right)^{q} - \left(\mu_{ikt}^{s}\right)^{q}\right)$$

Also since

$$|\Phi| = \left|\lambda_{t}\sum_{t=1}^{m}\left(\left(\mu_{ijt}^{s}\right)^{q} + \left(\mu_{jkt}^{s}\right)^{q} + \left(v_{ikt}^{s}\right)^{q} - \left(v_{ijt}^{s}\right)^{q} - \left(v_{jkt}^{s}\right)^{q} - \left(\mu_{ikt}^{s}\right)^{q}\right)\right|$$

$$\leq \lambda_{t}\sum_{t=1}^{m}\left|\left(\mu_{ijt}^{s}\right)^{q} + \left(\mu_{jkt}^{s}\right)^{q} + \left(v_{ikt}^{s}\right)^{q} - \left(v_{ijt}^{s}\right)^{q} - \left(v_{jkt}^{s}\right)^{q} - \left(\mu_{ikt}^{s}\right)^{q}\right|$$

And

$$CI(A) = \frac{1}{l}\frac{2}{n(n-1)(n-2)}\sum_{1\leq i<j<k}^{n}\sum_{s=1}^{l}|\Phi|,$$

So we can get

$$CI(A) \leq \frac{1}{l}\frac{2}{n(n-1)(n-2)}\sum_{1\leq i<j<k}^{n}\sum_{s=1}^{l}\left(\lambda_{t}\sum_{t=1}^{m}\left|\left(\mu_{ijt}^{s}\right)^{q} + \left(\mu_{jkt}^{s}\right)^{q} + \left(v_{ikt}^{s}\right)^{q} - \left(v_{ijt}^{s}\right)^{q} - \left(v_{jkt}^{s}\right)^{q} - \left(\mu_{ikt}^{s}\right)^{q}\right|\right)$$

$$= \lambda_{t}\sum_{t=1}^{m}CI(A) \leq \lambda_{t}\sum_{t=1}^{m}\overline{CI} = \overline{CI}$$

Thus, according to **Definition 10**, the collective q-ROHFPR $A$ is additive consistent, which completes the proof.

In order to measure the consensus degree between individual evaluation opinions and group opinions, and obtain scientific and effective program ranking, based on the **Definition 11**, the distance of $A_t = (a_{ijt})_{n\times n}$ with $a_{ijt} = (\mu_{ijt}, v_{ijt})$ and $A = (a_{ij})_{n\times n}$ with $a_{ij} = (\mu_{ij}, v_{ij})$ can be shown in formula(17).

$$D(A_t, A) = \frac{1}{2l}\frac{1}{(n-1)(n-2)}\sum_{i,j=1,i<j}^{n}\sum_{s=1}^{l}\left\{\left|u_{ijt}^{(s)} - u_{ij}^{(s)}\right| + \left|v_{ijt}^{(s)} - v_{ij}^{(s)}\right|\right\}, \tag{17}$$

In formula (17), $t = 1, 2, \ldots, m$, $s = 1, 2, \ldots, l$, $\mu_{ijt} = \{\mu_{ijt}^{1}, \mu_{ijt}^{2}, \ldots, \mu_{ijt}^{l}\}$, $v_{ijt} = \{v_{ijt}^{1}, v_{ijt}^{2}, \ldots, v_{ijt}^{l}\}$, $\mu_{ij} = \{\mu_{ij}^{1}, \mu_{ij}^{2}, \ldots, \mu_{ij}^{l}\}$ and $v_{ij} = \{v_{ij}^{1}, v_{ij}^{2}, \ldots, v_{ij}^{l}\}$.

According to formula (17), the distance-based consensus index is defined as follows.

**Definition 15**: Let $A_t = (a_{ijt})_{n \times n}$ ($t = 1, 2, \ldots, m$) be q-ROHFPRs with acceptable additive consistency, $A = (a_{ij})_{n \times n}$ be the aggregation matrix of $A_t$ obtained by **Definition 15**, and the consensus index ($GCI$) of $A_t$ is shown in formula (18).

$$GCI(A_t) = D(A_t, A) = \frac{1}{2l} \frac{1}{(n-1)(n-2)} \sum_{i,j=1, i<j}^{n} \sum_{s=1}^{l} |a_{ijt} - a_{ij}| \qquad (18)$$

Substitute formula (21) into the above formula to obtain:

$$GCI(A_t) = \frac{1}{2l} \frac{1}{(n-1)(n-2)} \sum_{i,j=1, i<j}^{n} \sum_{s=1}^{l} \{|\mu_{ijt}^s - \mu_{ij}^s| + |v_{ijt}^s - v_{ij}^s|\} \qquad (19)$$

In formula (23), $t = 1, 2, \ldots, m$. $D(A_t, A)$ is the distance between $A_t$ and $A$. The smaller $GCI(A_t)$ is, the more consistent the opinions of the kth expert are with those of group decision-making. If $GCI(A_t) = 0$, then the opinions of the t-th expert are completely consistent with those of group decision-making. However, in practical decision-making, because it is difficult for experts to reach a complete consensus in q-ROHFPRs, experts should also give the threshold $\overline{GCI}$ of consensus index to judge whether q-ROHFPRs is an acceptable consensus.

**Definition 16**: If $GCI(A_t) \leq \overline{GCI}$, then $A_t$ is an acceptable consensus q-ROHFPR. Otherwise, $A_t$ is called q-ROHFPR with unacceptable consensus. Particularly, when $\overline{GCI} = 0$ was established, $A_t$ was consensus completely. Where $\overline{GCI}$ is the consistency threshold and $\overline{GCI} \in [0,1]$.

### 5.2 Iteration Algorithm of Consensus

The reaching consensus process is the key part of getting reasonable decisions in GDM, it is need to consider not only the additive consistency of each expert preference matrix, but also the consensus problem between a single preference matrix and the aggregation matrix. If the evaluation opinions given by decision-makers differ greatly, the calculated results will not be convincing. we will focus on the issue of consensus, Improve the consensus of decision information among experts and groups, and then put forward specific group decision-making steps.

In a GDM problem, if $GCI(A_t) > \overline{GCI}$, then $A_t$ is an unacceptable consensus matrix. At this point, $A_t$ must be returned to the decision maker $d_t$ to reconstruct a new $q - \text{ROHFPR}$ that is closer to the group q-rung orthopair hesitant fuzzy preference relation $(Gq-ROHFPR)$. This process will be repeated until the predefined consensus level is achieved, we call this process the consensus-reaching process.

In the following, an automatic iterative algorithm is developed to describe the consensus-reaching process. Let $\tilde{A}_t = (\tilde{a}_{ijt})_{n\times n}$ be an acceptable consensus q-ROHFPR, take $A$ as the reference for optimizing consensus Algorithm.

**Input**: Individual q-ROHFPRs $\tilde{A}_t = (\tilde{a}_{ijt})_{n\times n}$ $(t = 1,2,\ldots,m)$, the weight vector of DMs $\lambda = (\lambda_1, \lambda_2, \ldots, \lambda_m)^T$, the predefined threshold $\overline{GCI}$, the maximum number of iterative times $\theta_{max} \geq 1$ and $0 < \zeta < 1$.

**Output**: The number of iterations $\theta$, modified q-ROHFPRs $\tilde{A}_t^\theta = (\tilde{a}_{ijt}^\theta)_{n\times n}$ $(t = 1,2,\ldots,m)$, the group consensus q-ROHFPR $A^\theta = (a_{ij}^\theta)_{n\times n}$, the consensus index of each q-ROHFPR $GCI(\tilde{A}_t^\theta)$.

**Step 1**: Set $t = 0$ and $\tilde{A}_t^{(0)} = (\tilde{a}_{ijt}^{(0)})_{n\times n} = \bar{A} = (\bar{a}_{ijt})_{n\times n}$.

**Step 2**: Utilize Eq. (17) to fuse all of the individual q-ROHFPRs $\tilde{A}_t = (\tilde{a}_{ijt})_{n\times n}$ $(t = 1,2,\ldots,m)$ into a $Gq-ROHFPR$ $A = (a_{ij})_{n\times n}$. For convenience, let $A^{(0)} = (a_{ij}^{(0)})_{n\times n} = A = (a_{ij})_{n\times n}$, $\theta = 0$, and let $\overline{GCI}$ be the dead line of acceptable consensus between each individual decision matrix and the group decision matrix.

**Step 3**: Calculate the group consensus index $GCI(\tilde{A}_{t,(0)})(t = 1,2,\ldots,m)$. If $GCI(\tilde{A}_{t,(0)}) \leq \overline{GCI}$ for any $t = 1,2,\ldots,m$, or $\theta \geq \theta_{max}$, then go to Step 5; otherwise, go to Step 4.

**Step 4**: Let $A_t^{(\theta+1)} = (a_{ijt}^{(\theta+1)})_{n\times n}$ and $A^{(\theta+1)} = (a_{ij}^{(\theta+1)})_{n\times n}$, where

$$a_{ijt}^{(\theta+1)} = \zeta a_{ijt}^{(\theta)} + (1-\zeta)a_{ij}^{(\theta)}, \quad t = 1,2,\ldots,m, i,j = 1,2,\ldots,n,$$

If $CI\left(A_t^{(\theta+1)}\right) \leq \overline{CI}$, set $\theta = \theta + 1$ and go to Step 2. Otherwise, use the model $(M-2)$ to adjust the consistency, then go to step 2.

**Step 5**: Output the number of iterations $\theta$, modified q-ROHFPRs $\tilde{A}_t = (\tilde{a}_{ijt})_{n\times n}$ $(t = 1,2,\ldots,m)$, the group consensus q-ROHFPR $A^\theta = (a_{ij}^\theta)_{n \times n}$, and the consensus index of each q-ROHFPR $GCI(\tilde{A}_t^\theta)$.

**Step 6:** End.

The above iterative algorithm is convergent. In fact, we can prove the following theorem:

**Theorem 2:** Let $A_t = (a_{ijt})_{n \times n}$ $(t = 1, 2, \ldots, m)$ be the individual q-ROHFPRs, let $A_t^\theta$ be a q-ROHFPR sequence that is generated by Algorithm 3 for the DM $d_t$, $GCI(A_t^\theta)$ be the consensus index of $A_t^\theta$ then for each $\theta$,

$$GCI\left(A_t^{(\theta+1)}\right) = \zeta GCI\left(A_t^{(\theta)}\right) < GCI\left(A_t^{(\theta)}\right), \ 0 < \zeta < 1, \ t = 1, 2, \ldots, m$$

and especially, let $\alpha^* = 0$, then $\lim\limits_{\theta \to \infty} GCI(A_t^{(\theta)}) = 0$, $t = 1, 2, \ldots, m$.

Proof. For $GCI\left(A_t^{(\theta)}\right) > \alpha^*$, according to Definition **16**,

$$GCI\left(A_t^{(\theta+1)}\right) = \frac{1}{2l} \frac{1}{(n-1)(n-2)} \sum_{1 \leq i < j}^{n} \sum_{s=1}^{l} \left| a_{ijt}^{(\theta+1)} - a_{ij}^{(\theta+1)} \right|$$

$$= \frac{1}{2l} \frac{1}{(n-1)(n-2)} \sum_{1 \leq i < j}^{n} \sum_{s=1}^{l} \left| a_{ijt}^{(\theta+1)} - \sum_{t=1}^{m} \lambda_t a_{ijt'}^{(\theta+1)} \right|$$

$$= \frac{1}{2l} \frac{1}{(n-1)(n-2)} \sum_{1 \leq i < j}^{n} \sum_{s=1}^{l} \left| \sum_{t'=1}^{m} \lambda_{t'} \left( a_{ijt}^{(\theta+1)} - a_{ijt'}^{(\theta+1)} \right) \right|$$

$$= \frac{1}{2l} \frac{1}{(n-1)(n-2)} \sum_{1 \leq i < j}^{n} \sum_{s=1}^{l} \left| \sum_{t'=1}^{m} \lambda_{t'} \left[ \left( \zeta a_{ijt}^{(\theta)} + (1-\zeta) a_{ij}^{(\theta)} \right) - \left( \zeta a_{ijt'}^{(\theta)} + (1-\zeta) a_{ij}^{(\theta)} \right) \right] \right|$$

$$= \frac{1}{2l} \frac{1}{(n-1)(n-2)} \sum_{1 \leq i < j}^{n} \sum_{s=1}^{l} \left| \sum_{t'=1}^{m} \lambda_{t'} \left( \zeta a_{ijt}^{(\theta)} - \zeta a_{ijt'}^{(\theta)} \right) \right|$$

$$= \frac{\zeta}{2l} \frac{1}{(n-1)(n-2)} \sum_{1 \leq i < j}^{n} \sum_{s=1}^{l} \left| \sum_{t'=1}^{m} \lambda_{t'} \left( a_{ijt}^{(\theta)} - a_{ijt'}^{(\theta)} \right) \right|$$

$$= \frac{\zeta}{2l} \frac{1}{(n-1)(n-2)} \sum_{1 \leq i < j}^{n} \sum_{s=1}^{l} \left| a_{ijt}^{(\theta)} - \sum_{t'=1}^{m} \lambda_{t'} a_{ijt'}^{(\theta)} \right|$$

$$= \frac{\zeta}{2l} \frac{1}{(n-1)(n-2)} \sum_{1 \leq i < j}^{n} \sum_{s=1}^{l} \left| a_{ijt}^{(\theta)} - a_{ij}^{(\theta)} \right|$$

$$= \zeta GCI(A_t^{(\theta)}) < GCI(A_t^{(\theta)})$$

i.e., $GCI\left(A_t^{(\theta+1)}\right) < GCI\left(A_t^{(\theta)}\right)$, for all $t = 1, 2, \ldots, m$.

Afterward, since

$$GCI\left(A_t^{(\theta+1)}\right) = \zeta GCI\left(A_t^{(\theta)}\right) = \zeta^2 GCI\left(A_t^{(\theta-1)}\right) = \ldots = \zeta^\theta GCI\left(A_t^{(\theta)}\right) = \zeta^{\theta+1} GCI(A_t), \text{ for all } t = 1, 2, \ldots, m.$$

And if $\alpha^* = 0$, we get

$$\lim_{\theta \to \infty} GCI(A_t^{(\theta+1)}) = \lim_{\theta \to \infty} \zeta^{\theta+1} GCI(A_t) = GCI(A_t) \lim_{\theta \to \infty} \zeta^{\theta+1} = 0, \text{ for all } t = 1, 2, \dots, m.$$

which completes the proof.

Obviously, the outstanding feature of the algorithm is that it can automatically modify different individual opinions in order to reach a consensus in the group opinions, and avoid forcing experts to modify their decision matrix, so as to provide great convenience for its practical application.

### 5.3 GDM Algorithm of q-ROHFPRs

This section describes the MOGDM problems under q-rung orthopair hesitant fuzzy environments. Using the above analysis, a new method is presented to solve such GDM problems with q-ROHFPRs. The concrete steps of the proposed method for GDM with q-ROHFPRs are described in the following.

**Step 1.** Suppose $\overline{CI}$ is the threshold of the acceptable additive consistency. Calculate the consistency index of each q-ROHFPR $A_t = (a_{ijt})_{n \times n}$ $(t = 1, 2, \dots, m)$ using Eq.(9). According to **Definition 11**, check whether $A_t$ satisfies acceptable additive consistency. If all of them are acceptably consistent, go to Step 4. Otherwise, go to Step 3.

**Step 2:** Use the model $(M - 2)$ to adjust each q-ROHFPR $A_t$ to satisfies the acceptable additive consistency threshold, and then the adjusted matrix is expressed as $\tilde{A}_t$, $t = (1, 2, \dots, m)$.

**Step 3:** Aggregate the q-ROHFPRs $\tilde{A}_t = (\tilde{a}_{ijt})_{n \times n}$ $(t = 1, 2, \dots, m)$ obtained in Step 3 into the collective q-ROHFPR $A = (a_{ij})_{n \times n}$ based on Eq.(16).

**Step 4:** Let $\overline{GCI}$ be the consensus threshold, where $\overline{GCI} \in [0,1]$. The value of the group consensus index $GCI(\tilde{A}_t)$ are generated utilizing Eq. (19). As described in **Definition 16**, if $GCI(\tilde{A}_t) \leq \overline{GCI}$ means $\tilde{A}_t$ has acceptable consensus, set $\tilde{A}_{t'} = \tilde{A}_t$ and go to Step 7. Otherwise, skip to the Step 6.

**Step 5:** Use the iteration algorithm in Section 5.2, the decision group is improved to the acceptable group consensus and the adjusted q-ROHFPRs $\tilde{A}_t = (\tilde{a}_{ijt})_{n \times n}$ $(t = 1, 2, \dots, m)$ are obtained.

**Step 6:** As per Eq.(16), the final aggregation matrix $A_G$ of the adjusted q-ROHFPRs $\tilde{A}_t$ is determined. Utilize model $(M-4)$ to derive the q-rung orthopair hesitant fuzzy weight vector $w = (w_1, w_2, ..., w_t)^T$ of the alternatives.

**Step 7:** Calculate the score function $S(w_t)$ ($t = 1, 2, ..., m$) by Eq.(2).

**Step 8:** Rank the score value $S(w_t)$ and get the ranking order of alternatives. The larger the value of $S(w_t)$, the better the alternative.

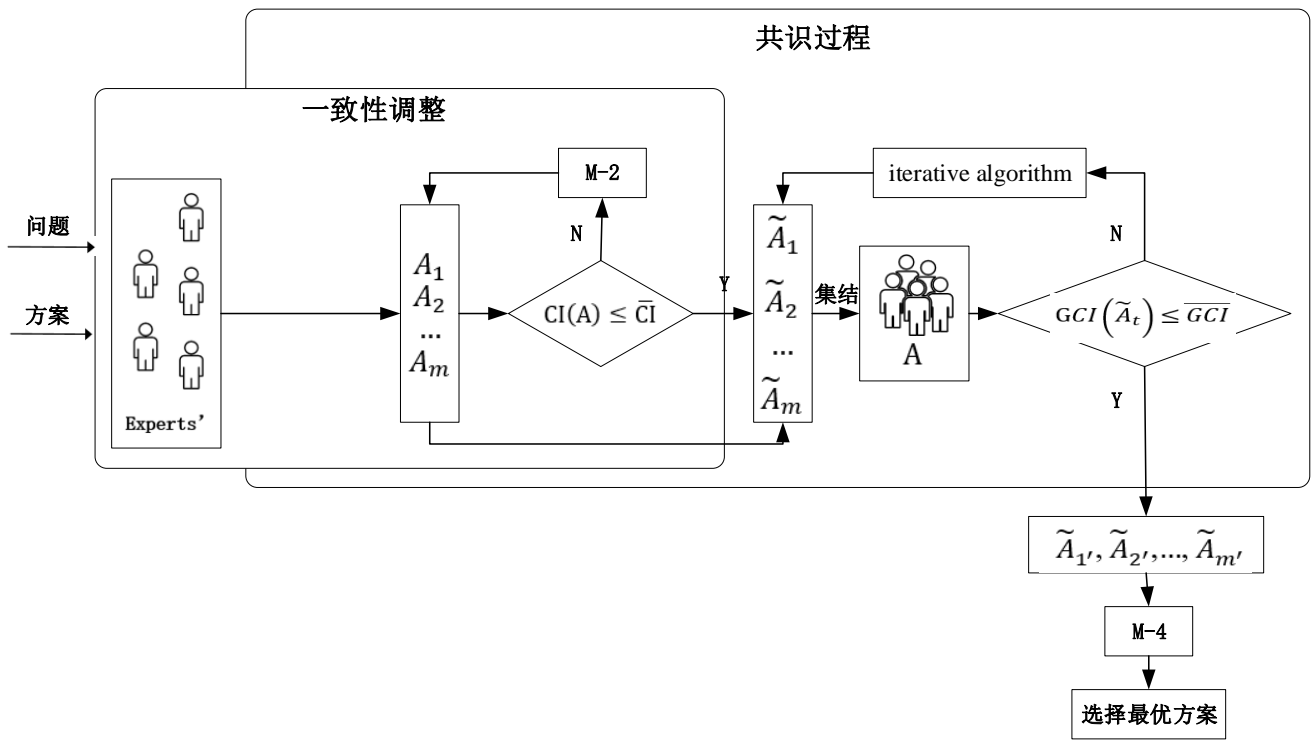

Fig.2 Decision making process for GDM with q-ROHFPRs

# 6 Numerical example

## 6.1 Numerical example

Due to the frequent failure of the laboratory computer equipment, the School of Software and Internet of Things Engineering decided to purchase a batch of new equipment to replace it. There are currently 4 computer brands (alternatives) $X = \{x_1, x_2, x_3, x_4\}$ to choose from. The college invites 3 experts $D = \{d_1, d_2, d_3\}$ to evaluate the alternatives from the four attributes of hardware, budget, computer after-sales, and brand influence, and each expert's evaluation information Both are represented by $q-ROHFPRs$ $A^k (k = 1,2,3)$, the expert weight is $\lambda = (0.3, 0.2, 0.5)$, because the cost and benefit type do not require complex calculations and the amount of data is small, this paper directly gives the converted evaluation matrix.

$A_1$

$$= \begin{pmatrix} <\{\sqrt[q]{0.5}\},\{\sqrt[q]{0.5}\}> & <\{0.5333,0.6000,0.6667\}, \{0.1167,0.1500,0.1833\}> & <\{0.3333,0.4001,0.4667\}, \{0.4167,0.4500,0.4833\}> & <\{0.3166,0.3499,0.3832\}, \{0.2500,0.3500,0.4500\}> \\ <\{0.1167,0.1500,0.1833\}, \{0.5333,0.6000,0.6667\}> & <\{\sqrt[q]{0.5}\},\{\sqrt[q]{0.5}\}> & <\{0.2334,0.3001,0.3668\}, \{0.4167,0.4500,0.4833\}> & <\{0.6167,0.6500,0.6833\}, \{0.1333,0.2001,0.2667\}> \\ <\{0.4167,0.4500,0.4833\}, \{0.3333,0.4000,0.4667\}> & <\{0.4167,0.4500,0.4833\}, \{0.2334,0.3001,0.3668\}> & <\{\sqrt[q]{0.5}\},\{\sqrt[q]{0.5}\}> & <\{0.3333,0.4001,0.4667\}, \{0.3333,0.4001,0.4667\}> \\ <\{0.2500,0.3500,0.4500\}, \{0.3166,0.3499,0.3833\}> & <\{0.1333,0.2001,0.2667\}, \{0.6167,0.6500,0.6833\}> & <\{0.3333,0.4001,0.4667\} \{0.3333,0.4001,0.4667\}> & <\{\sqrt[q]{0.5}\},\{\sqrt[q]{0.5}\}> \end{pmatrix}$$

$A_2$

$$= \begin{pmatrix} <\{\sqrt[q]{0.5}\},\{\sqrt[q]{0.5}\}> & <\{0.2334,0.3001,0.3668\}, \{0.4167,0.4500,0.4833\}> & <\{0.3166,0.3499,0.3832\}, \{0.4167,0.4500,0.4833\}> & <\{0.2167,0.2500,0.2833\}, \{0.3166,0.3499,0.3832\}> \\ <\{0.4167,0.4500,0.4833\}, \{0.2334,0.3001,0.3668\}> & <\{\sqrt[q]{0.5}\},\{\sqrt[q]{0.5}\}> & <\{0.7167,0.7500,0.7833\}, \{0.1167,0.1500,0.1833\}> & <\{0.3500,0.4500,0.5500\}, \{0.2334,0.3001,0.3668\}> \\ <\{0.4167,0.4500,0.4833\}, \{0.3166,0.3499,0.3832\}> & <\{0.1167,0.1500,0.1833\}, \{0.7167,0.7500,0.7833\}> & <\{\sqrt[q]{0.5}\},\{\sqrt[q]{0.5}\}> & <\{0.4167,0.4500,0.4833\}, \{0.1500,0.2500,0.3500\},> \\ <\{0.3166,0.3499,0.3832\}, \{0.2167,0.2500,0.2833\}> & <\{0.2334,0.3001,0.3668\}, \{0.3500,0.4500,0.5500\}> & <\{0.1500,0.2500,0.3500\}, \{0.4167,0.4500,0.4833\}> & <\{\sqrt[q]{0.5}\},\{\sqrt[q]{0.5}\}> \end{pmatrix}$$

$A_3$

$$= \begin{pmatrix} <\{\sqrt[q]{0.5}\},\{\sqrt[q]{0.5}\}> & <\{0.3333,0.4001,0.4667\}, \{0.2334,0.3001,0.3668\},> & <\{0.4167,0.4500,0.4833\}, \{0.3166,0.3499,0.3832\}> & <\{0.4167,0.4500,0.4833\}, \{0.2167,0.2500,0.2833\}> \\ <\{0.2334,0.3001,0.3668\}, \{0.3333,0.4001,0.4667\}> & <\{\sqrt[q]{0.5}\},\{\sqrt[q]{0.5}\}> & <\{0.5167,0.5500,0.5833\}, \{0.1500,0.2500,0.3500\}> & <\{0.7167,0.7500,0.7833\}, \{0.1167,0.1500,0.1833\}> \\ <\{0.3166,0.3499,0.3832\}, \{0.4167,0.4500,0.4833\}> & <\{0.1500,0.2500,0.3500\}, \{0.5167,0.5500,0.5833\}> & <\{\sqrt[q]{0.5}\},\{\sqrt[q]{0.5}\}> & <\{0.4334,0.5000,0.5667\}, \{0.2167,0.2500,0.2833\}> \\ <\{0.2167,0.2500,0.2833\}, \{0.4167,0.4500,0.4833\}> & <\{0.1167,0.1500,0.1833\}, \{0.7167,0.7500,0.7833\}> & <\{0.2167,0.2500,0.2833\}, \{0.4334,0.5000,0.5667\}> & <\{\sqrt[q]{0.5}\},\{\sqrt[q]{0.5}\}> \end{pmatrix}$$

**Step 1:** Check the consistency of the q-ROHFPRs $A_t$, and calculate the consistency index $CI(A_t)$ of q-ROHFPRs to get: $CI(A_1) = 0.1611$, $CI(A_2) = 0.125$, $CI(A_3) = 0.0833$. It can be seen from $CI(A_1) > 0.1$, $CI(A_2) > 0.1$, $CI(A_3) < 0.1$ that $A^1$ and $A^2$ do not meet the consistency conditions, but $A^3$ is consistent.

**Step 2:** Use $(M-2)$ to improve the consistency of $A^1$, $A^2$, and get $\tilde{A}^1$ and $\tilde{A}^2$.

$\tilde{A}_2$

$$= \begin{pmatrix} <\{\sqrt[q]{0.5}\},\{\sqrt[q]{0.5}\}> & <\{0.5333,0.6000,0.6667\}, \{0.1167,0.1500,0.1833\}> & <\{0.3333,0.4001,0.4667\}, \{0.4167,0.4500,0.4833\}> & <\{0.3166,0.3499,0.3832\}, \{0.2500,0.3500,0.4500\}> \\ <\{0.1167,0.1500,0.1833\}, \{0.5333,0.6000,0.6667\}> & <\{\sqrt[q]{0.5}\},\{\sqrt[q]{0.5}\}> & <\{0.2334,0.3001,0.3668\}, \{0.4167,0.4500,0.4833\}> & <\{0.6167,0.6500,0.6833\}, \{0.1333,0.2001,0.2667\}> \\ <\{0.4167,0.4500,0.4833\}, \{0.3333,0.4001,0.4667\}> & <\{0.4167,0.4500,0.4833\}, \{0.2334,0.3001,0.3668\}> & <\{\sqrt[q]{0.5}\},\{\sqrt[q]{0.5}\}> & <\{0.3333,0.4001,0.4667\}, \{0.3333,0.4001,0.4667\}> \\ <\{0.2500,0.3500,0.4500\}, \{0.3166,0.3499,0.3832\}> & <\{0.1333,0.2001,0.2667\}, \{0.6167,0.6500,0.6833\}> & <\{0.3333,0.4001,0.4667\} \{0.3333,0.4001,0.4667\}> & <\{\sqrt[q]{0.5}\},\{\sqrt[q]{0.5}\}> \end{pmatrix}$$

$\tilde{A}_2$

$$= \begin{pmatrix} <\{\sqrt[q]{0.5}\},\{\sqrt[q]{0.5}\}> & <\{0.2334,0.3001,0.3668\}, \{0.4167,0.4500,0.4833\}> & <\{0.3166,0.3499,0.3832\}, \{0.4167,0.4500,0.4833\}> & <\{0.2167,0.2500,0.2833\}, \{0.3167,0.3499,0.3832\}> \\ <\{0.4167,0.4500,0.4833\}, \{0.2334,0.3001,0.3668\}> & <\{\sqrt[q]{0.5}\},\{\sqrt[q]{0.5}\}> & <\{0.7167,0.7500,0.7833\}, \{0.1167,0.1500,0.1833\}> & <\{0.3500,0.4500,0.5500\}, \{0.2334,0.3001,0.3668\}> \\ <\{0.4167,0.4500,0.4833\}, \{0.3166,0.3499,0.3832\}> & <\{0.1167,0.1500,0.1833\}, \{0.7167,0.7500,0.7833\}> & <\{\sqrt[q]{0.5}\},\{\sqrt[q]{0.5}\}> & <\{0.4167,0.4500,0.4833\}, \{0.1500,0.2500,0.3500\},> \\ <\{0.3167,0.3499,0.3832\}, \{0.2167,0.2500,0.2833\}> & <\{0.2334,0.3001,0.3668\}, \{0.3500,0.4500,0.5500\}> & <\{0.1500,0.2500,0.3500\}, \{0.4167,0.4500,0.4833\}> & <\{\sqrt[q]{0.5}\},\{\sqrt[q]{0.5}\}> \end{pmatrix}$$

**Step 3:** Use formula (16) to calculate the aggregation matrix $A$ of $\tilde{A}_1$, $\tilde{A}_2$ and $A_3$, after calculation, $A$ is obtained as follows:

$$A = \begin{pmatrix} <\{\sqrt[q]{0.5}\},\{\sqrt[q]{0.5}\}> & <\{0.3522,0.4217,0.4904\}, \{0.2930,0.3328,0.3760\},> & <\{0.3616,0.4062,0.4500\}, \{0.3810,0.4139,0.4469\}> & <\{0.3224,0.3568,0.3910\}, \{0.2598,0.3154,0.3766\}> \\ <\{0.2930,0.3328,0.3760\}, \{0.3522,0.4217,0.4904\}> & <\{\sqrt[q]{0.5}\},\{\sqrt[q]{0.5}\}> & <\{0.4382,0.4928,0.5443\}, \{0.2930,0.3330,0.3832\}> & <\{0.5673,0.6278,0.6844\}, \{0.1680,0.2226,0.2784\}> \\ <\{0.3810,0.4139,0.4469\}, \{0.3616,0.4062,0.4500\}> & <\{0.2930,0.3330,0.3832\}, \{0.4382,0.4928,0.5443\}> & <\{\sqrt[q]{0.5}\},\{\sqrt[q]{0.5}\}> & <\{0.3950,0.4535,0.5115\}, \{0.2569,0.3145,0.3766\}> \\ <\{0.2598,0.3154,0.3766\}, \{0.3224,0.3568,0.3910\}> & <\{0.1680,0.2226,0.2784\}, \{0.5673,0.6278,0.6844\}> & <\{0.2569,0.3145,0.3766\}, \{0.3950,0.4535,0.5115\}> & <\{\sqrt[q]{0.5}\},\{\sqrt[q]{0.5}\}> \end{pmatrix}$$

**Step 4:** Consensus judgments are made on $\tilde{A}_1$, $\tilde{A}_2$ and $A_3$, and the consensus degree $GCI(A_t)$ between the q-ROHFPRs $A_t$ of each expert and the aggregation matrix $A$ is calculated respectively. Substitute $\tilde{A}_1$, $\tilde{A}_2$ and $A_3$, and $A$ into the formula (21) to obtain $GCI(\tilde{A}^1) = 0.08, GCI(\tilde{A}^2) = 0.1, GCI(A_3) = 0.06$. Since $GCI(\tilde{A}^1) < 0.1, GCI(\tilde{A}^2) = 0.1, GCI(A_3) < 0.1$ satisfy the consensus conditions, no adjustment is required.

**Step 5:** Calculate the weight $w = (w_1, w_2, ..., w_k)^T$ of the alternatives according to the aggregation matrix $A$, ($i = 1,2,...,n; s = 1,2,...,l$).

$$\min \sum_{i=1}^{n-1} \sum_{j=i+1}^{n} \sum_{s=1}^{l} \lambda_{ij}^{s+} + \lambda_{ij}^{s-} + \theta_{ij}^{s+} + \theta_{ij}^{s-}$$

$$s.t.\begin{cases} u_{ij}^s - \sqrt[q]{0.5(w_{iu}^s)^q + 0.5(w_{jv}^s)^q} - \lambda_{ij}^{s+} + \lambda_{ij}^{s-} = 0, \\ v_{ij}^s - \sqrt[q]{0.5(w_{iv}^s)^q + 0.5(w_{ju}^s)^q} - \theta_{ij}^{s+} + \theta_{ij}^{s-} = 0, \\ 0 \leq w_{iu}^s, w_{iv}^s \leq 1, \\ 0 \leq (w_{iu}^s)^q + (w_{iv}^s)^q \leq 1, \\ \sum_{j=1,j\neq i}^{n}(w_{iu}^s)^q \leq (w_{iv}^s)^q, \\ (w_{iu}^s)^q + n - 2 \geq \sum_{j=1,j\neq i}^{n}(w_{iv}^s)^q, \\ \lambda_{ij}^{s+}, \lambda_{ij}^{s-}, \theta_{ij}^{s+}, \theta_{ij}^{s-} \geq 0, \lambda_{ij}^{s+} \times \lambda_{ij}^{s-} = 0, \theta_{ij}^{s+} \times \theta_{ij}^{s-} = 0 \end{cases}$$

Get:

$$w_1 = (\{0.2178, 0.2916, 0.4667\}, \{0.2755, 0.3935, 0.4500\})^T$$
$$w_2 = (\{0.0006, 0.0016, 0.0007\}, \{0.4984, 0.4950, 0.5126\})^T$$
$$w_3 = (\{0.2755, 0.3935, 0.4500\}, \{0.3994, 0.4692, 0.4667\})^T$$
$$w_4 = (\{0.2178, 0.2916, 0.4500\}, \{0.3759, 0.4065, 0.4823\})^T$$

**Step 7:** Substitute $w = (w_1, w_2, w_3, w_4)^T$ into formula (2), we get:

$$S(w_1) = -0.253, \ S(w_2) = -0.6009, \ S(w_3) = -0.7613, \ S(w_4) = -1.0273$$

**Step 8:** Sort the alternatives. Since $S(w_1) > S(w_2) > S(w_3) > S(w_4)$, $x_1$ is optimal.

## 6.2 Comparative analysis

When the parameter α is unchanged and q takes different values, we get the corresponding calculation results and ranking order, as shown in Table 2. It can be seen from Table 2 that for different values of q, the ranking order of the alternatives may be the same, which means that no matter the risk attitude of DMs is risky or neutral, and the ranking order of the alternatives is always $x_1 > x_2 > x_3 > x_4$. The above results show that the changes of parameters q has little effect on the decision-making results, indicating that the algorithm is stable.

Table 2. When $α = 0.5$ ranking orders of alternatives for different values of parameter $q$

|  | $x_1$ | $x_2$ | $x_3$ | $x_4$ | 结果 |
| --- | --- | --- | --- | --- | --- |
| q=3 | -0.253 | -0.6009 | -0.7613 | -1.0273 | $x_1 > x_2 > x_3 > x_4$ |
| q=4 | -0.1001 | -2572 | -0.3142 | -0.4228 | $x_1 > x_2 > x_3 > x_4$ |

# 7 Conclusion

This paper proposes q-ROHFPRs to solve the GDM problems. Based on q-ROHFS, the definition of q-ROHFPR is proposed, and several properties of q-ROHFPR are studied in detail. The definitions of additive consistency and acceptable additive consistency of q-ROHFPRs are proposed, and a consistency optimization model is designed to improve the consistency of q-ROHFPRs. This model is calculated on the basis of the original q-ROHFPR itself. No other q-ROHFPR is involved, and the preference information conveyed by the original q-ROHFPR can be retained to the greatest extent. The q-rung orthopair hesitant fuzzy weighted average is used to aggregate multiple q-ROHFPRs into one Gq-ROHFPR. This operator can guarantee the consistency of Gq-ROHFPR. Through individual and group q-ROHFPRs, the acceptable consensus is defined, and a consensus improvement model is constructed to adjust the degree of consensus between individuals and groups. Use the group q-ROHFPR to generate the q-ROHF priority vector of the alternatives, and then use the q-ROHFN scoring function to calculate the scores of different computer brands, which is convenient for ranking and selecting the best. The group decision-making method proposed in this paper provides a new method for the laboratory to purchase computer brands. At the same time, compared with other decision-making methods, it is concluded that the decision-making method proposed in this paper is more stable and beneficial to brand evaluation. For this reason, the proposed group decision-making method provides a reliable guarantee for the laboratory's decision to purchase computer equipment.

The group decision problem studied in this paper involves relatively little research data. With the development of fuzzy theory, the application fields involved are more and more extensive, including text recognition, image and speech, natural language, earthquake prediction, etc., and the types of data are increasing. The more diverse, so the group decision problem that contains a large amount of data will be studied in the near future, and we will also apply the method proposed in this paper to other application areas, such as hospital management, pattern recognition and artificial intelligence.

## References

[1] Zadeh, L. A. Fuzzy sets[J]. Information and Control, 1965,8(3):338–353.
[2] Torra V., Hesitant fuzzy sets[J]. International Journal of Intelligent Systems, 2010,25:529-539.
[3] Saraji, M. K., Mardani, A., Köppen, M., Mishra, A. R., & Rani, P. An extended hesitant fuzzy set using SWARA-MULTIMOORA approach to adapt online education for the control of the pandemic spread of COVID-19 in higher education institutions[J]. Artificial Intelligence Review, 2021. https://doi.org/10.1007/s10462-021-10029-9
[4] Zindani, D., Maity, S. R., Bhowmik, S. Extended TODIM method based on normal wiggly hesitant fuzzy sets for deducing optimal reinforcement condition of agro-waste fibers for green product development[J]. Journal of Cleaner Production, 2021, 301:126947.
[5] Liu, X., Wang, Z., Zhang, S., & Garg, H. Novel correlation coefficient between hesitant fuzzy sets with application to medical diagnosis[J]. Expert Systems with Applications, 2021,183: 115393.
[6] Hao, Z., Xu, Z., Zhao, H., & Su, Z. Optimized data manipulation methods for intensive hesitant fuzzy set with applications to decision making[J]. Information Sciences, 2021,580:55–68.
[7] Li D, Zeng W, Li J. New distance and similarity measures on hesitant fuzzy sets and their applications in multiple criteria decision making[J]. Eng Appl Artif Intell, 2015,40:11-16.

[8] M. M. Xia and Z. S. Xu. Hesitant fuzzy information aggregation in decision making[J]. International Journal of Approximate Reasoning, 2011,52(3): 395–407.

[9] Z. S. Xu, M. M. Xia. Distance and similarity measures for hesitant fuzzy sets[J]. Information Sciences, 2011,181(11): 2128–2138.

[10] Zhu B, Xu Z, Xia M. Dual hesitant fuzzy sets[J]. J Appl Math,2012, 1–12.

[11] Chen, N., Xu, Z., Xia, M. Interval-valued hesitant preference relations and their applications to group decision making[J]. Knowledge-Based Systems, 2013, 37:528–540.

[12] D. Liu, D. Peng and Z. Liu. The distance measures between q-rung orthopair hesitant fuzzy sets and their application in multiple criteria decision making[J]. *Int. J. Intell. Syst.*, 2019, 34(9): 2104-2121.

[13] Saaty, TL. The Analytic Hierarchy Process McGraw-Hill[J].1980.

[14] Fanyong Meng,Shyi-Ming Chen. A framework for group decision making with multiplicative trapezoidal fuzzy preference relations[J]. Information Sciences,2021,577:722-747.

[15] Zhang, Z., Chen, S.-M. Group decision making based on multiplicative consistency-and-consensus preference analysis for incomplete q-rung orthopair fuzzy preference relations[J]. Information Sciences, 2021,574:653–673.

[16] Xu ZS, Liao HC. The multiplicative consistency index of hesitant fuzzy preference relation[J]. IEEE Transactions on Fuzzy Systems 2015, doi: 10.1109/TFUZZ.2015.2426315.

[17] Zhu, B. , Xu, Z. Regression methods for hesitant fuzzy preference relations[J]. Technological and Economic Development of Economy, 2013,19 (sup1), S214–S227 .

[18] Zhao, Na; Xu, Zeshui; Liu, Fengjun. Group Decision Making with Dual Hesitant Fuzzy Preference Relations[J]. Cognitive Computation, 2016,8(6), 1119–1143.

[19] Chen N, Xu ZS, Xia MM. Correlation coefficients of hesitant fuzzy sets and their applications to clustering analysis[J]. Appl Math Model ,2013,37:2197–2211

[20] Tang, Jie; Meng, Fanyong. New method for interval-valued hesitant fuzzy decision making based on preference relations. Soft Computing,2020,24(17): 13381-13399.

[21] Tang, J., Zhang, Y., Fujita, H., Zhang, X., Meng, F. Analysis of acceptable additive consistency and consensus of group decision making with interval-valued hesitant fuzzy preference relations[J]. Neural Computing and Applications, 2021,33(13), 7747–7772.

[22] Zhang Z , Kou X , Dong Q . Additive consistency analysis and improvement for hesitant fuzzy preference relations[J]. Expert Systems with Applications, 2018, 98:118-128.

[23] Zhu B, Xu ZS, Xu JP. Deriving a ranking from hesitant fuzzy preference relations under group decision making[J]. IEEE Trans Cybern, 2014,44(8):1328–1337.

[24] Li CC, Rodrı́guez RM, Herrera F, Martı́nez L, Dong YC. Consistency of hesitant fuzzy linguistic preference relations: an interval consistency index[J]. Inf Sci,2018,432:361–437.

[25] Tang J, An QX, Meng FY, Chen XH. A natural method for ranking objects from hesitant fuzzy preference relations[J]. Int J Inf Technol Decis Making, 2017, 16(6):1611–1646.

[26] Zhang Z, Kou XY, Dong QX. Additive consistency analysis and improvement for hesitant fuzzy preference relations[J]. Expert Syst Appl,2018, 98:118–128

[27] Zhang Z, Kou XY, Yu WY, Guo CH. On priority weights and consistency for incomplete hesitant fuzzy preference relations[J]. Knowl-Based Syst ,2018,141(1):115–126.

[28] Meng, F. Y., An, Q. X. A new approach for group decision making method with hesitant fuzzy preference relations[J]. Knowledge-Based Systems, 2017,127(Suppl. C), 1–15.

[29] Wu Z , Jin B , Fujita H. Consensus analysis for AHP multiplicative preference relations based on consistency control: A heuristic approach[J]. Knowledge-Based Systems, 2019, 191.

[30] Zhang Z , Pedrycz W . Iterative Algorithms to Manage the Consistency and Consensus for Group Decision-

Making with Hesitant Multiplicative Preference Relations[J]. IEEE Transactions on Fuzzy Systems, 2019, (99):1-1.

[31] He, Y., Xu, Z. S. A consensus reaching model for hesitant information with different preference structures[J]. Knowledge-Based Systems, 2017,135:99–112.

[32] Zhang C, Liao HC, Li L. Additive consistency‑based priority‑generating method of q‑rung orthopair fuzzy preference relation[J]. Int J Intell Syst. 2019;34(9):2151‑2176.

[33] Tanino, T.: Fuzzy preference orderings in group decision making[J]. Fuzzy Sets Syst. 1984,12:117–131.

[34] Zhang, Z., Chen, S. M. (2021). Group decision making with incomplete q-rung orthopair fuzzy preference relations[J]. Information Sciences, 553, 376–396.